\documentclass[twoside,a4paper,12pt,centertags,reqno]{amsart} 
\usepackage{amsmath,amssymb,verbatim,vmargin}
\usepackage{color}
\usepackage{tikz}
\usepackage{hyperref}
\hypersetup{colorlinks,bookmarks=true,linktocpage=true,citecolor=blue,linkcolor=magenta}

\usepackage{eulervm}   

\allowdisplaybreaks 
\usepackage{color}
\usepackage[all]{xy} 

\usepackage{mathtools}
\usepackage{MnSymbol} 

\theoremstyle{plain}
\newtheorem{thm}{Theorem}[section]
\newtheorem{lem}[thm]{Lemma}

\theoremstyle{definition}

\newtheorem{rem}[thm]{Remark}

\usepackage{enumerate}

\newcommand{\ep}{\varepsilon} 

\newcommand{\bC}{{\mathbb C}}

\newcommand{\bN}{{\mathbb N}}

\def\barint_#1{\mathchoice
            {\mathop{\vrule width 6pt
height 3 pt depth -2.5pt
                    \kern -9.5pt
\intop \kern -4pt}\nolimits_{#1}}%
            {\mathop{\vrule width 5pt height
3 pt depth -2.6pt
                    \kern -6.5pt
\intop \kern -4pt}\nolimits_{#1}}%
            {\mathop{\vrule width 5pt height
3 pt depth -2.6pt
                    \kern -6pt
\intop \kern -4pt}\nolimits_{#1}}%
            {\mathop{\vrule width 5pt height
3 pt depth -2.6pt
          \kern -6pt \intop \kern -4pt}\nolimits_{#1}}}
          
           \def\bariint_#1{\mathchoice
            {\mathop{\vrule width 15pt
height 3 pt depth -2.5pt
                    \kern -15.8pt
\intop \kern -8pt\intop \kern -4pt}\nolimits_{#1}}%
            {\mathop{\vrule width 9pt height
3 pt depth -2.6pt
                    \kern -10.5pt
\intop \kern -8pt\intop \kern -4pt}\nolimits_{#1}}%
            {\mathop{\vrule width 9pt height
3 pt depth -2.6pt
                    \kern -10pt
\intop \kern -8pt\intop \kern -4pt}\nolimits_{#1}}%
            {\mathop{\vrule width 9pt height
3 pt depth -2.6pt
          \kern -8pt \intop \kern -10pt\intop \kern -4pt}
      \nolimits_{  #1}}}

\def\barintlim_#1{\mathchoice
            {\mathop{\vrule width 6pt
height 3 pt depth -2.5pt
                    \kern -8.8pt
\intop \kern -4pt}\limits_{#1}}%
            {\mathop{\vrule width 5pt height
3 pt depth -2.6pt
                    \kern -6.5pt
\intop \kern -4pt}\limits_{#1}}%
            {\mathop{\vrule width 5pt height
3 pt depth -2.6pt
                    \kern -6pt
\intop \kern -4pt}\limits_{#1}}%
            {\mathop{\vrule width 5pt height
3 pt depth -2.6pt
          \kern -6pt \intop \kern -4pt}\limits_{#1}}}
          
           \def\bariintlim_#1{\mathchoice
            {\mathop{\vrule width 15pt
height 3 pt depth -2.5pt
                    \kern -15.8pt
\intop \kern -8pt\intop \kern -4pt}\limits_{#1}}%
            {\mathop{\vrule width 9pt height
3 pt depth -2.6pt
                    \kern -10.5pt
\intop \kern -8pt\intop \kern -4pt}\limits_{#1}}%
            {\mathop{\vrule width 9pt height
3 pt depth -2.6pt
                    \kern -10pt
\intop \kern -8pt\intop \kern -4pt}\limits_{#1}}%
            {\mathop{\vrule width 9pt height
3 pt depth -2.6pt
          \kern -8pt \intop \kern -10pt\intop \kern -4pt}
      \limits_{  #1}}}
          
\renewcommand{\iint}{\int \kern -3pt\int}       

\numberwithin{equation}{section}
\usepackage[symbol]{footmisc}

\title{On Galbis' integration lemmas}

\author{Yi C. Huang} 
\address{School of Mathematical Sciences, Nanjing Normal University, Nanjing 210023, People's Republic of China}
\email{Yi.Huang.Analysis@gmail.com}
\urladdr{https://orcid.org/0000-0002-1297-7674}

\author{Fei Xue} 
\address{School of Mathematical Sciences, Nanjing Normal University, Nanjing 210023, People's Republic of China}
\email{05429@njnu.edu.cn}

\date{\today} 

\thanks{Research of the authors is partially supported by the National NSF grants of China (nos. 11801274 and 12201307) and the Jiangsu Provincial NSF grant (no. BK20210555).
YCH thanks Jian-Yang Zhang for switching \textit{counter-clockwise} the seminar material from Nicola-Tilli to Galbis.}

\begin{document}

\begin{abstract}
We simplify in this note Galbis' proof of certain norm estimates for self-adjoint Toeplitz operators on the Fock space.
This relies on an extension (and a unification) of his integration lemmas, yet with a simpler proof in the same spirit.
\end{abstract}

\maketitle


\section{Introduction}

Assuming that the bounded symbol is further radial and integrable,
Galbis \cite{Gal22} obtained some very interesting norm estimates for self-adjoint Toeplitz operators on the Fock space on $\bC$.
See also Grudsky and Vasilevski \cite{GruVas02} for a related result.
Galbis' arguments rely crucially on the following two elementary integration lemmas.

\begin{lem} \label{lem1}
Let $I\subset[0,\infty)$ be a measurable set with finite Lebesgue measure. Then
$$\frac{1}{n!}\int_I s^ne^{-s}ds\leq 1-e^{-|I|}.$$
\end{lem}

\begin{lem} \label{lem2}
Let $(I_k)_{k=1}^N$ be disjoint sets with finite measure and $0\leq\ep_k\leq1$ for every $1\leq k\leq N$. 
Then for every $p\in\bN_0=\{0, 1, 2, \cdots\}$ we have
$$\sum_{k=1}^N\ep_k\int_{I_k} \frac{t^p}{p!}e^{-t}dt\leq 1-\exp\left(-\sum_{k=1}^N\ep_k|I_k|\right).$$
\end{lem}

The aim of this note is to point out the following extension of Lemmas \ref{lem1}-\ref{lem2}.

\begin{lem} \label{lem3}
Let $d\mu(s)=g(s)ds$ with $g\geq0$ integrable and $\|g\|_\infty\leq1$. Then
$$\sup_{n\in\bN_0}\frac{1}{n!}\int_0^\infty s^ne^{-s}d\mu(s)\leq 1-e^{-\|g\|_{L^1(0,\infty)}}.$$
\end{lem}

\begin{rem}
Galbis' arguments for \cite[Theorem 1]{Gal22} are now greatly simplified by our Lemma \ref{lem3}:
using his notations, for $|F(z)|=g(|z|)$ and $\sum_{p=0}^\infty |b_p|^2=1$, we have
$$\begin{aligned}
\sum_{p=0}^\infty |b_p|^2\int_0^\infty g\left(\sqrt{t/\pi}\right)\frac{t^p}{p!}e^{-t}dt
&\leq \sup_{p\in\bN_0}\int_0^\infty g\left(\sqrt{t/\pi}\right)\frac{t^p}{p!}e^{-t}dt\\
&\leq 1-e^{-\left\|g\left(\sqrt{\cdot/\pi}\right)\right\|_{L^1(0,\infty)}}=1-\exp\left(-\|F\|_{L^1(\bC)}\right).
\end{aligned}$$
Here $F$ is the Toeplitz symbol and $F=\sum_{p=0}^\infty b_pe_p$, where $e_p(z)=(\pi^p/p!)^{1/2}z^p$.
Moreover, this approach also enables us to bypass an approximation argument on $g$ (in connection with a result by Hu and Lv \cite{HuLv14} for the Toeplitz operators).
\end{rem}

We apologise to the reader for the necessary briefness of this short note and suggest he (or she) has (at least) Galbis' article \cite{Gal22} handy.

\section{Proof of Lemma \ref{lem3}}

The proof is in the same spirit of Galbis' proof of Lemma \ref{lem1},
and is reminiscent of Hardy's integration lemma \cite[Proposition 3.6, page 56]{BenSha88}.  
Indeed, given $n\in\bN_0$, $h_n(s):=\frac{s^n}{n!}e^{-s}$ attains its absolute maximum at $s=n$.
Moreover, $h_n$ increases on $[0,n]$ and decreases on $[n,\infty)$.
Note also that $g\geq0$ and $\|g\|_\infty\leq1$.
So the maximum of the integral $\frac{1}{n!}\int_0^\infty s^ne^{-s}g(s)ds$ is attained while
$g$ is the indicator function of some interval $[a,b]$\footnotemark\footnotetext{Note that $[a,b]$ depends on $n$.} 
that contains $n$ and has length $\|g\|_{L^1(0,\infty)}$. 
Therefore,   
$$\frac{1}{n!}\int_0^\infty s^ne^{-s}d\mu(s)\leq \int_a^bh_n(s)ds\leq 1-e^{-\|g\|_{L^1(0,\infty)}}.$$
In the second inequality we use Galbis' nice estimation in his proof of Lemma \ref{lem1}:
$$\begin{aligned}
\int_a^bh_n(s)ds&=\frac{e^{-a}}{n!}\int_0^{b-a}(t+a)^ne^{-t}dt\\
&=\sum_{k=0}^nC_n^k\frac{a^{n-k}}{n!}e^{-a}\int_0^{b-a}t^ke^{-t}dt\\
&=\sum_{k=0}^n\frac{a^{n-k}}{(n-k)!}e^{-a}\frac{1}{k!}\int_0^{b-a}t^ke^{-t}dt\\
&\leq \sup_{0\leq k\leq n}\frac{1}{k!}\int_0^{b-a}t^ke^{-t}dt\\
&=\sup_{0\leq k\leq n}\left(1-e^{-(b-a)}\sum_{j=0}^k\frac{(b-a)^j}{j!}\right)=1-e^{-(b-a)}.
\end{aligned}$$
The lemma is then proved by varying $n\in\bN_0$.

\bigskip

\section*{\textbf{Compliance with ethical standards}}

\bigskip

\textbf{Conflict of interest} The authors have no known competing financial interests
or personal relationships that could have appeared to influence this reported work.

\bigskip

\textbf{Availability of data and material} Not applicable.

\bigskip

\bibliographystyle{alpha}

\bibliography{HuaY-XueF-GalbisLemma} 

\begin{thebibliography}{Gal22}

\bibitem[BS88]{BenSha88}
Colin Bennett and Robert~C. Sharpley.
\newblock {\em Interpolation of Operators}.
\newblock Academic Press, 1988.

\bibitem[Gal22]{Gal22}
Antonio Galbis.
\newblock Norm estimates for selfadjoint {T}oeplitz operators on the {F}ock
  space.
\newblock {\em Complex Analysis and Operator Theory}, 16(1):15, 2022.

\bibitem[GV02]{GruVas02}
Sergei~M. Grudsky and Nicolai~L. Vasilevski.
\newblock Toeplitz operators on the {F}ock space: radial component effects.
\newblock {\em Integral Equations and Operator Theory}, 44(1):10--37, 2002.

\bibitem[HL14]{HuLv14}
Zhangjian Hu and Xiaofen Lv.
\newblock Toeplitz operators on {F}ock spaces {$F^p(\varphi)$}.
\newblock {\em Integral Equations and Operator Theory}, 80(1):33--59, 2014.

\end{thebibliography}
 
\end{document}